\documentclass[11pt]{article}
\usepackage{amsmath,amssymb}
\usepackage{hyperref}
\usepackage{jkss,here}
\usepackage{natbib}
\submit{final} \volumn{30}{2} 

\def\anticlockwise{{\circlearrowleft}}
\def\clockwise{{\circlearrowright}}
\def\ket#1{{ | #1 \rangle }}
\def\bra#1{{ \langle #1 | }}
\def\braket#1#2{{ \langle #1  |#2 \rangle }}
\def\ketbra#1#2{{\ket{#1}\bra{#2}}}
\def\smile{{\mbox{\tt :\!-)\;}}}
\def\frown{{\mbox{\tt :\!-(\;}}}

\newtheorem{thm}{Theorem}

\begin{document}
    
\heading{Teleportation into Quantum Statistics}{Richard D. Gill}

\title{Teleportation into Quantum Statistics}

\author{Richard Gill\footnote{Mathematical Institute, Leiden University; gill@math.leidenuniv.nl;
\url{https://www.math.leidenuniv.nl/~gill}. The paper was written while I was at Utrecht University.}}
\begin{abstract}
The paper is a tutorial introduction to quantum information theory, 
developing the basic model and emphasizing the role of statistics and 
probability.
\end{abstract}

\section*{PRELUDE}

Between the present prelude and a concluding postlude, 
the body of this paper is divided into five numbered sections. 

For motivation and introduction, Section 1 contains a discussion of recent experiments 
in solid state physics, constructing a single bit (0/1 memory register) 
of a new kind of computer called the quantum computer. 

In Section 2 we give the mathematical model behind quantum 
computation, quantum communication, quantum statistics, quantum probability, 
quantum stochastics; the whole field now being called \emph{quantum information theory}.
We will see that the model is (mathematically speaking) elementary, 
it is essentially probabilistic, and it leads to natural statistical 
problems. The model is built on precisely four ingredients: notions of 
(i) \emph{state} of a quantum system, (ii) its time 
\emph{evolution}, (iii) the formation of joint systems 
from separate, called \emph{entanglement}, and finally, (iv), the stochastic 
interface with the real world: \emph{measurement}. In Section 2
we restrict attention to basic forms of these notions: 
states are actually so-called pure states, (represented by vectors); 
evolutions are unitary; measurements are  so-called simple measurements 
(projector-valued probability measures). 
Later we will see how combining these building blocks in various ways leads to 
generalized notions of state, evolution and measurement. But the four 
ingredients in their basic form remain the only items on which the whole theory 
is built.

In Section 3, an intermission, we will illustrate the basic model ingredients
with the example of quantum teleportation, 
which in a few lines of elementary 
algebra and a simple probability calculation exemplifies all the key 
model ingredients, the statistical challenges, and the extraordinary
physical implications, of the theory.

In Section 4 we take a new look at the 
model ingredients, extending the notions of pure states and simple measurements 
to mixed states (density matrices) and generalized measurements 
(operator-valued probability measures; more generally, completely 
positive instruments). This enables us to describe the problems of
quantum statistical design and quantum statistical inference in a compact 
and precise way, and it also gives hope that these problems might have elegant 
solutions.

In Section 5 we develop some theory of quantum statistical inference.
For a quantum statistical model, we define
the quantum score and the quantum Fisher information, 
leading to quantum Cram\'er-Rao and information bounds (by now, very 
classical material). We briefly survey some recent progress in  
quantum statistical design and inference, in particular the quantum 
information bound of \citet{gillmassar00} and their results on 
asymptotically optimal quantum design and inference. This gives 
solutions to problems posed by the motivating example of Section 1: 
how can the experimentalists substantiate their claims, with a 
minimum of experimental effort?

In a postlude or maybe more appropriately, aftermath, 
we will switch to a more polemical mode and comment on the 
relations between quantum probability, quantum statistics, quantum physics and 
technology, and real probability and real statistics. 

The aim is to convince the reader that the area of quantum statistical 
inference is grounded in a simple mathematical model, which combines 
basic elements from probability, statistics, linear algebra, and (the 
absolutely basic) elements of complex analysis. A bit of trigonometry 
also comes in handy. No physics knowledge at all, is needed.
Most statisticians' training includes all of these ingredients. 
However many will not have been exposed to what comes out of the ``intersections''
between these fields, for instance, linear algebra with vectors and 
matrices of complex numbers instead of real numbers.
But one just needs to learn a few useful facts about eigenvalues and 
eigenvectors of complex matrices, 
which directly generalize the familiar facts about symmetric real 
matrices to \emph{self-adjoint} complex matrices, and generalize 
real orthonormal matrices to complex \emph{unitary} matrices.

One does not need any physics 
background to appreciate the basic modelling, and from there, to 
contribute to the scientific development of quantum information theory.
The field is associated with some of the most significant 
current developments in physics, of deep scientific importance and holding promise 
of substantial technological impact. The physics we are talking about has an 
essential probabilistic component; the experiments which are being done now and 
the experiments which will play a role in developing the new technologies, are going 
to need statistical design and analysis. 
Starting from the basic model, one can quickly pose intriguing problems 
of statistical design and inference, some of which have elegant and exciting 
solutions, others are quite open. These toy problems 
are related to current work in physics, information theory, and computer science, 
at this moment of great theoretical interest, and likely to become of practical 
interest in the near future. Already, people working in computer science 
and in information theory have turned in a big way to (theoretical) 
quantum computation and quantum information theory. For instance, in 
Korea I refer to the work of Dong Pyo Chi and his colleagues at Seoul 
National University. Strangely, probabilists 
and statisticians do not seem to be making a similar move. We will give some 
thoughts on why this should be so, in the postlude. 

The survey papers \citet{gill01} and \citet{barndorffnielsenetal01} 
cover many further topics, especially drawing attention to open problems, 
and moreover give further references, especially for background reading.

\section{MOTIVATING EXAMPLE: THE DELFT QUBIT}\label{sec:1}

We briefly discuss an experiment carried out in Delft, the Netherlands, by the 
group of Prof.\ Hans Mooij. See {\tt http://vortex.tn.tudelft.nl/},
and especially the pictures on the personal pages of Ph.D student Caspar van 
der Wal. The experiment was reported in {\it Science\/} in 1999. 
Switching on a magnetic field causes electric current 
to flow around a superconducting aluminium ring. The aluminium ring is a 
thousandth of a millimeter in diameter, and a billion electrons are involved in the 
current flow.  From a classical physical viewpoint one can imagine just two kinds 
of current flow of a given size in this little circuit: {\it 
clockwise}, $\clockwise$,
and {\it anti-clockwise}, $\anticlockwise$. 
The claim of the experimenters was that they produced an electric current in the 
state 
$\ket\psi~=~\alpha\ket\clockwise+\beta\ket\anticlockwise$, 
where $\alpha$ and $\beta$ are two complex numbers, with $|\alpha|^2 +|\beta|^2=1$.    
$\ket\clockwise$ and  $\ket\anticlockwise$ stand for two orthogonal 
unit vectors in a complex vector space, we can think of them as two-dimensional 
complex column vectors of length $1$, say the basis vectors, 
$
\begin{pmatrix}
1\\ 
0
\end{pmatrix}
$ and 
$
\begin{pmatrix}
0\\ 
1
\end{pmatrix}
$. 
This object has been called {\it The Delft Qubit\/}; a qubit being a single bit in the memory of 
a future quantum computer. A classical computer works with a memory, the bits 
of which can register only $0$ or $1$, however a quantum computer allows coherent 
superpositions of $0$ and $1$, such as the state I have just talked about. Another 
description is {\it The Schr\"odinger Squid\/}; this name refers to the device: a 
Superconducting Quantum Interference Device; and to the infamous Schr\"odinger 
cat.  Now one might ask, how could the experimenters know that this state has 
been produced? Well, by repeating the experiment about ten thousand times, and 
each time measuring the current. This is done by a second squid surrounding the 
first, and connected to the outside world by a lot of circuitry. It does not directly 
give us estimates of $\alpha$ and $\beta$.  In fact, in first instance, it does nothing 
interesting at all: the measurement essentially looks to see whether the current is 
flowing $\clockwise$ or $\anticlockwise$. This forces the quantum state to jump into 
either of the states $\ket\clockwise$ or $\ket\anticlockwise$, 
and it makes this choice with 
probabilities $|\alpha|^2$ and $|\beta|^2$. The experimentalists find the same 
values of these 
probabilities (relative frequencies), as are predicted by an elaborate theoretical 
physical calculation concerning the whole system.

So this does not {\it prove} anything at all: one would have seen the same relative 
frequencies, if the qubit had from the start been, in a fraction $|\alpha|^2$ of the times, 
in state $\ket\clockwise$, and in a fraction $|\beta|^2$ of the times, 
in state $\ket\anticlockwise$.  
However, small developments in the technology of this experiment will make the 
finding more secure. The aim is not just to create qubits but to manipulate them. 
In particular, it should be possible to implement the 
transformation of the state, which sends the original orthonormal basis vectors  
$
\begin{pmatrix}
1\\ 
0
\end{pmatrix}
$ 
and 
$
\begin{pmatrix}
0\\ 
1
\end{pmatrix}
$, 
into the new orthonormal basis vectors
$
\begin{pmatrix}
1/\sqrt{2}\\ 
1/\sqrt{2}
\end{pmatrix}
$ 
and 
$
\begin{pmatrix}
1/\sqrt{2}\\ 
-1/\sqrt{2}
\end{pmatrix}
$. 
The result of this {\it unitary transformation\/} is to convert the original qubit into the 
state $\frac 1{\sqrt 2}(\alpha+\beta)\ket\clockwise+
\frac 1{\sqrt 2}(\alpha-\beta)\ket\anticlockwise$. If we now 
measure, we will find relative frequencies of $|\alpha+\beta|^2 /2$ 
and $|\alpha-\beta|^2$, 
different from the relative frequencies had the state been initially in a fraction 
$|\alpha|^2$ of the times, $|\clockwise\rangle$, and in a fraction $|\beta|^2$  of the times, 
$|\anticlockwise\rangle$. (As the reader may compute, one would then 
have observed  $\clockwise : \anticlockwise$ in equal proportions).

Still, the experiment is difficult to do. The question considered in 
this paper is: what is the optimal experimental design in order to 
determine the actually created state of this quantum qubit as 
accurately as possible with as small as possible number of 
repetitions of the experiment? The answer to this question was still 
completely unknown two years ago; in fact, the correct answer turns 
out to be quite opposite to many physicists' intuition as to what is 
best. But the question needs also to be further specified, since what is 
best depends on what experimental resources are available, what 
prior knowledge there is about the state being measured, and the 
relative importance of different features of the state.

Building a quantum qubit is just a first step towards building a 
quantum computer. Though many technologies are being explored for 
this purpose (ion traps, nuclear magnetic resonance, optics) the Delft 
implementation has the promise of scalability: the possibility to 
control not one or two but thousands of qubits.
The idea of quantum computation is to store program and data of some 
algorithm, coded in $0$'s and $1$'s, into the states $\ket 0$ and 
$\ket 1$ of a number of qubits. 
The whole system then evolves unitarily, and at the end of this evolution, a series 
of (possibly random) $0$'s and $1$'s are read off by measuring each qubit 
separately. The possibilities allowed by the basic model of quantum mechanics 
allow, for instance, (with an algorithm of Peter Shor) to factor large integers in 
polynomial time, which will make all currently used cryptography methods 
obsolete! Fortunately quantum cryptography promises a secure alternative. One 
cannot look at a qubit without disturbing it, and if this idea is cleverly exploited, it 
becomes possible to transmit messages coded in qubit states in such a way, that 
the interference of any eavesdropper would be detected by the recipient. 
Quantum computation may still be far away, and moreover it is not entirely 
clear if it would live up to its promises. But there is a 
strong feeling that quantum optical communication technology is just around 
the corner. In any case, as integrated circuits become smaller and 
communication speeds faster, present-day technology is rapidly approaching 
quantum limits. On the other hand, new quantum technologies can exploit precisely 
those phenomena that for the older technologies is a barrier to further progress.

\section{THE BASIC MODEL INGREDIENTS}

What is the basic mathematical model behind all this, what then are the statistical 
problems, and what do we know about the solutions? We have seen the notion of 
{\it states\/} (more precisely, pure state), mathematically formalized 
as vectors $\ket\psi$ in a complex vector space, of unit length: 
$\braket\psi\psi=1$. States can be {\it unitarily transformed}, that is to say, one 
may implement an orthonormal transformation (change of basis) and get a new 
state. In  principle, any desired unitary transformation could be implemented by 
setting up appropriate external fields. It is a manipulation of the state of the 
quantum system, involving, for instance, magnetic fields, which one can control, 
but without back-action on the real world outside. No information passes from 
the quantum system into the real world. What we have not yet described is the 
mathematical model for bringing initially separate quantum systems into 
(potential) interaction with one another. This is the essential ingredient of the 
quantum computer: one should not have $N$ separate qubits, but one quantum 
system of $N$ qubits together. The appropriate model for this is the 
{\it formation of tensor products}. In words, two separate systems 
brought together have as state, a 
vector in a space of dimension equal to the product of the two original 
dimensions; and the new state vector has as components, all the products 
of two
components, one from each of the two original state vectors. The $N$ qubits of a quantum 
computer live in a $2^N$ dimensional state space. The initial state is a product state, 
but a unitary evolution can bring the joint system into a state, which cannot be 
represented as a product state. This phenomenon is called entanglement.

The last ingredient has already been touched upon, and that is {\it measurement}. At 
this stage, and only at this stage, is information passed from the quantum system 
into the real world. The information is random, and its probability distribution 
depends on the state of the system. The system makes a random jump to a new 
state. The basic measurement is characterized by a collection of orthogonal 
subspaces of the state space, together spanning the whole space; and a real 
number or label, associated to each subspace. The collection of subspaces and 
numbers corresponds to an experiment one might do in the laboratory. When the 
experiment is carried out, the state vector of the quantum system is projected 
into one of the subspaces (and renormalised to have length one); the 
corresponding label becomes known in the real world; and all this happens with 
probability equal to the squared length of the projection of the original state 
vector into the subspace. By Pythagoras, these squared lengths add up to 1.

These are all the ingredients: state vectors (also called pure states), unitary 
evolution, entanglement (formation of product systems), and (simple) 
measurement. We now go through them more formally, giving as special 
example the important case of a two-dimensional state space:
this applies to the qubit, to a two-level system, to polarization of a 
photon, to spin of an electron or other `spin-half' particles.

\subsection{States}

The only definition of a quantum system is: a physical system which 
satisfies the laws of quantum mechanics, and those are the laws which 
we are about to outline. According to modern physics, quantum 
mechanics rules at all levels: atoms, molecules; light, 
electromagnetic radiation; the 
early universe (cosmology); string theorists apply it to fundamental 
constituents of matter at much lower scale (much higher energy level) 
than anything which is nowadays attainable by experiment. In any case, 
it is a physical system (or certain aspects of a physical system) 
whose interaction with the rest of the world is so simple that it can 
be succesfully described according to the following picture. The state 
of the system is: precisely what you need to know, in order to make 
predictions about the results of any future experiments with the 
system. These predictions are probabilistic, so to be more precise: 
when we know the state of a system, we know the relative frequency of 
the possible outcomes of any possible measurement on the system, in 
many repetitions of the experiment: do such and such a measurement on 
identical copies of a system is such and such a state. Identical 
copies just means: prepared in identical fashion.

In this section we will 
represent the state of a quantum system with a non-zero complex vector.
In all our examples, the state space will be finite 
dimensional, say $d$-dimensional, and a state vector is therefore just a column
vector of $d$ complex numbers. (More generally one needs to work in a 
separable Hilbert space). We will use both notations $\psi$ and 
$\ket\psi$ to stand for the state vector. The adjoint of this 
vector is the row vector containing the complex conjugates of the 
elements of $\psi$. It is denoted by $\psi^{*}$ or by $\bra\psi$, 
again two 
notations for precisely the same thing. It follows that $\psi^{*}\psi$, or if 
you prefer $\braket\psi\psi$, stands for the squared length of the 
vector $\psi$ (the sum of squared abolute values of its elements).
If $\psi$ is a state-vector, then all the non-zero vectors in the 
one-dimensional subspace $[\psi]=\{z\psi:z\in\mathbb C\}$ actually 
represent the same state (i.e, the physical predictions are identical).
Conventionally, one normalizes state vectors to have length $1$, thus
$\braket\psi\psi=1$. It is then easy to check that the matrix 
$\rho=\psi\psi^{*}=\ketbra\psi\psi=\Pi_{[\psi]}$ is the matrix which 
orthogonally projects a arbitrary vector to the subspace $[\psi]$.
Since one can reconstruct the subspace $[\psi]$ from knowing the 
matrix $\rho$, it follows that one can equally well represent states 
by the matrix $\rho$ as by the vector $\psi$. Even if $\psi$ is 
normalized, one can still multiply the state-vector by a complex 
number of absolute value $1$, i.e., of the form $e^{i\theta}$ for some 
real angle $\theta\in[0,2\pi)$, and get a different vector, which is 
also a representative of the same state. The angle $\theta$ is called 
a phase. So an overall phase is irrelevant, but when writing one state 
vector as a linear combination of others, the relative phases do make 
a difference.

Note how the at first sight clumsy notation $\ket\psi$, $\bra\psi$, 
helps one to graphically recognise whether one is talking about a 
number $\braket\psi\psi$ or a matrix $\ketbra\psi\psi$. A further 
advantage is that we are now 
also able to denote state vectors by replacing the name of a 
vector, $\psi$, with a verbal or graphic description of the 
state, as in $\ket{\mathrm{Alive}}$ and $\ket{\mathrm{Dead}}$, or 
$\ket{\,\smile\,}$ and $\ket{\,\frown\,}$. The notation 
is due to Dirac; $\ket\psi$, $\bra\psi$ are called a ket and a bra 
respectively.

\subsubsection{Example of states: the qubit} 

The same mathematical model of a two-dimensional quantum system 
applies to all kinds of physical systems: the current in the Delft 
qubit, the gound state versus first excited state of an atom at very 
low temperature, the spin of an electron or other so-called spin-half 
particle, the polarization of a photon.
Whatever the application,
the state vector of a two dimensional 
quantum system can be written as $\alpha\ket 0+\beta\ket 1$ where $\ket 
0$ and $\ket 1$ are a pair of orthonormal basis 
elements of $\mathbb C ^{2}$, and $\alpha$ and $\beta$ are two 
complex numbers, not both zero. The labels $0$ and $1$ are 
conventionally used when talking about the quantum qubit (a single quantum 
memory bit). In other contexts other descriptive labels might be 
appropriate, as we have seen above. Normalizing the length of the 
vector to $1$, and taking the coefficient of $\ket 0$ to be a real 
number (which can be achieved by a suitable phase factor) one easily 
sees that one can represent the state by the vector $\cos\theta\ket 0
+\sin\theta e^{i\phi}\ket 1$, for some real angles $\theta\in[0,\pi]$ and 
$\phi\in[0,2\pi)$. We will see in a moment, that it is very useful to 
think of the angles $(\theta,\phi)$ as polar coordinates of a real
three-dimensional unit vector $\vec u$: $\theta$ is the co-latitude, i.e., 
the angle you have to move down from the North pole, $\phi$ is the 
longitude, the angle you have to move around the globe. When we are 
talking about spin of an electron (`spin-half') the direction of 
$\vec u$ in real three-space really can be thought of, as the 
direction of the axis of spin of the electron. In other applications 
(e.g., polarization of a photon, see Section 3)
the interpretation might be more complicated. But the mathematics is 
the same.
To know the state, one should equivalently specify a complex 
$2$-vector $\ket\psi$, real polar coordinates $(\theta,\phi)$,
or a real unit $3$-vector $\vec u$.
One might denote the state vector correspondingly as $\ket\psi$,
as $\ket{\theta,\phi}$ or as $\ket{\vec u}$. In the important 
application of a spin-half particle, e.g., an electron, the basis 
states are denoted $\ket\uparrow$ and $\ket\downarrow$, `up' and `down' 
respectively, and the state $\ket{\vec u}$ really can be thought of as
the state of an electron spinning in the real spatial direction $\vec 
u$.

The matrix representation of the same state is 
found by some simple algebra and trigonometry to be equal to $\rho(\theta,\phi)
=\rho(\vec u) = \frac12(\mathbf 1 +\vec u(\theta,\phi)\cdot\vec{\sigma}))$, 
where the ingredients in this formula are described as follows.
Bold type indicates complex two by two matrices which otherwise might be
confused with numbers. Thus
$\mathbf 1$ is the two by two identity matrix. The arrow indicates a 
vector of $3$ components, which might be reals or matrices.
$\vec u(\theta,\phi)$ is the real three-dimensional unit vector having
polar coordinates $(\theta,\phi)$. The symbol `$\cdot$' denotes the ordinary inner product, and
$\vec{\sigma}$ denotes a vector of three two by two matrices,
the so-called Pauli spin matrices, 
$\sigma_{x}=
\begin{pmatrix}
    0 & 1 \\ 
    1 & 0
\end{pmatrix}$,
$\sigma_{y}=
\begin{pmatrix}
    0 & -i \\ 
    i & 0
\end{pmatrix}$,
$\sigma_{z}=
\begin{pmatrix}
    1 & 0 \\ 
    0 & -1
\end{pmatrix}$.
So writing $\vec u=(u_{x},u_{y},u_{z})$, by definition $\vec u 
\cdot\vec{\sigma}=u_{x}\sigma_{x}+u_{y}\sigma_{y}+u_{z}\sigma_{z}$.
Each of the Pauli spin matrices is self-adjoint, $\sigma=\sigma^{*}$, where 
the adjoint of a matrix is the transpose of the matrix of complex conjugates
of the original matrix elements. Self-adjoint complex matrices, like real 
symmetric matrices, have real eigenvalues and an orthonormal basis of 
eigenvectors. In particular, the Pauli spin matrices all have eigenvalues
$+1$ and $-1$, their 
eigenvectors are $\psi(\pm\vec u_{x})$, $\psi(\pm\vec u_{y})$,
and $\psi(\pm\vec u_{z})$, where $\vec u_{x}$ denotes the real 
three-dimensional unit vector in the $x$-direction, and so on.
The \emph{opposite} real three-vectors $\vec u$ and $-\vec u$
correspond to \emph{orthogonal} state vectors $\ket{\vec u}$, $\ket{-\vec u}$.
Some useful properties of the spin matrices are $\sigma_{x}\sigma_{y}=
-\sigma_{y}\sigma_{x}=i\sigma_{z}$ (and the same for cyclic permuations of 
$(x,y,z)$), and $\sigma_{x}^{2}=\sigma_{y}^{2}=\sigma_{z}^{2}=\mathbf 1$.

Later we will extend from the so-called pure states, represented by a 
state vector $\psi$, to the what are called mixed states: the state 
vector of the quantum system is drawn with probability distribution 
$\mathrm P(\mathrm d \psi)$ from the set of all state vectors, let us
suppose them to be all normalized to length $1$.
It turns out (as we will see in Section 4) that for all physical predictions, 
it suffices to know no more and no less than the ordinary probability mixture 
$\rho_{\mathrm{ave}}=\int\rho(\psi)\mathrm P(\mathrm d 
\psi)$ of the corresponding state-matrices 
$\rho(\psi)=\ketbra\psi\psi$.
This simple mathematical fact has an extraordinary consequence. 
Suppose I give you a stream of spin-half particles, each independently 
prepared with equal probability in the state $\ket{\vec u_{z}}$ or in the state 
$\ket{-\vec u_{z}}$ (`up' and `down'). Or, I give you a stream of spin-half particles, 
each independently prepared with equal probability in the state 
$\ket{\vec u_{x}}$ or in the state 
$\ket{-\vec u_{x}}$ (`left' and `right'). Later under the subsection 
on
measurement, we will see how such a preparation could be made. The mixed 
state matrix for the first case is $\frac12(\rho(\vec u_{z})+\rho(-\vec 
u_{z}))$, for the second case it is $\frac12(\rho(\vec u_{x})+\rho(-\vec 
u_{x}))$, in both cases this average state matrix is the rather simple 
$\frac12\mathbf1$. Thus whatever measurements you make on the 
particles, you will never be able to tell the difference between the 
two scenarios. The statistical predictions of any experiment you can 
do, would be the same. This extraordinary fact casts doubt on the 
idea that the state of a quantum system, as some collection of real or 
complex numbers, is somehow `engraved' permanently on individual 
particles (electrons, photons, or whatever). If that were the case, it 
would be very strange that one could never decide whether a huge 
number of particles, each engraved with very different states, could 
never be distinguished. It seems that the state is not a property of 
an individual particle, but rather of the way a particle is created, 
and carries merely statistical information. This fact bothers 
physicists, who are not fond of randomness, a lot,
but probabilists and statisticians should  
find it relatively easy to live with.

\subsection{Evolutions}

A quantum system not acting in any way on the 
external world, may be influenced by it, in the following way. For any 
particular situation the physicist will be able to write down a 
self-adjoint matrix $H$ called the Hamiltonian, or energy, and then 
the state at time $t$ of a quantum system is derived from the state at 
time $0$ by solving the differential equation $ i\hbar\mathrm d \psi /\mathrm 
d t = H \psi$. Here $\hbar$ is Planck's constant, a rather small 
quantity of work = energy times time, and the equation we have just 
written down is the famous Schr\"odinger equation. The point is, that the 
experimentalist might be able to arrange for the same quantum system 
to evolve under different Hamiltonians $H$, for instance if we are 
talking about spin of electrons, by appropriately setting 
up different external magnetic fields. If we are not 
talking about spin and magnetism, but about polarization of photons, 
passing light through various crystals might implement different 
Hamiltonian evolutions. 

For our finite dimensional 
quantum systems we can solve the equation explicitly as 
$\psi(t)=e^{Ht/i\hbar}\psi(0)$. Even more explicitly, one can write 
the matrix $H$ in terms of its eigenvalue-eigenvector decomposition as
$H=\sum_{a}a\ketbra a a$, where $a$ runs through the eigenvalues of 
$H$,
which are real, and $\ket{H=a}=\ket a$ is a convenient notation for: 
the normalized eigenvector corresponding to eigenvalue $a$. One should 
actually say: a normalized eigenvector, there is still an arbitrary 
phase factor. And now, since $e^{Ht/i\hbar}=\sum_{a}e^{
at/i\hbar}\ketbra a a$ one can solve the time evolution 
as $\ket{\psi(t)}=\sum_{a} e^{at/i\hbar} 
\braket a {\psi(0)}\, \ket a$. This shows that a given state can be 
expressed as a complex superposition of energy eigenstates. Each 
eigenstate on its own evolves in a rather boring way: according to the 
phase factor $e^{at/i\hbar}$. However linear combinations of 
eigenstates can express fascinating interference effects, as the 
relative phases of the component eigenstates change in time.

Now the matrix $U=U_{t}=e^{Ht/i\hbar}$ has the special property of 
being unitary: that means precisely that $UU^{*}=U^{*}U=\mathbf 1$. In 
other words, the transformation $\psi\mapsto U\psi$ is nothing more 
nor less than a change of (orthonormal) basis of our state-space. The 
key point for the applications is that any unitary matrix $U$ whatsoever 
is of the form $U=U_{t}=e^{Ht/i\hbar}$ for some Hamiltonian $H$ 
and time length $t$. Thus in principle, if one could implement any 
particular Hamiltonian in the laboratory, one can implement any 
unitary transformation of a state.

\subsubsection{Example of evolution: the qubit} 

The matrices 
$
\begin{pmatrix}
    1 & 0 \\ 
    0 & -1
\end{pmatrix}
$
and 
$
\begin{pmatrix}
    0 & 1 \\ 
    1 & 0
\end{pmatrix}
$
are both unitary, and therefore correspond to transformations of a 
quantum state that one might implement in a laboratory. The first maps
an arbitrary state vector $\alpha\ket0+\beta\ket 1$ into $\alpha\ket0-\beta\ket 
1$, a sign change,
the second maps $\alpha\ket0+\beta\ket 1$ into $\alpha\ket 1 +\beta\ket 
0$, the so-called spin-flip.

There is a beautiful connection between the unitary transformations of 
states in $\mathbb C^{2}$ and the orthorgonal rotations of 
corresponding unit vectors in $\mathbb R^{3}$, involving the Pauli 
spin matrices, but we do not need it here.

\subsection{Entanglement}

In ordinary probability theory there is a natural way to model the 
bigger probability space formed by performing independently two other 
probability experiments. There is a very analogous, and 
mathematically very natural operation, for modelling the bringing 
together of two independent and completely separate quantum systems into 
(potential) interaction with one another. The mathematical tool for 
this is the notion of tensor product. A quantum system with state 
vector $\psi$ in a $d$-dimensional state space, together with another 
system with state vector $\phi$ in a $d'$ dimensional state space, 
together form a quantum system in a $d\times d'$ dimensional state 
space, with state vector $\psi\otimes\phi$, by which we mean the 
vector containing each product of one of the $d$ elements
of the vector $\psi$ 
with one of the $d'$ elements of the vector $\phi$, 
arranged in some fixed order which 
suits you. Now this particular state is not very interesting: as we 
will see in the next subsection, when one does simultaneous measurements 
on each of the two subsystems, the outcomes are independent and 
distributed exactly as they would have been, considered entirely 
separately. But the point is that this boring product state is the 
state of the joint system, only at the precise moment when the two 
subsystems are brought together. From that moment they will evolve 
together under some Hamiltonian. And if that Hamiltonian is not of 
the boring form $H\otimes\mathbf 1' + \mathbf 1\otimes H'$  
(use your imagination to define the tensor product of matrices now, 
rather than of vectors), the joint state will evolve in some period 
of time into a new state in the huge $d\times d'$ dimensional space, 
with a state vector which cannot be written in the simple product form 
which it had at time $0$. 

Every state vector in the big product space can be written as a 
complex linear combination of product states. Whenever one needs 
a linear combination of more than one product, we call the joint state 
entangled. As we will see, such states have remarkable properties.

\subsubsection{Example of entanglement: the qubit} 

The quantum computer will be built 
of a large collection, say $N$, of simple two-level systems or qubits.
Thus the state of the whole system is a vector in a $2^{N}$ 
dimensional state space, including states which are not of the 
special form: each qubit separately in its own state. The idea of the 
quantum computer is to implement the logical transformations on bits, 
which are the basis of classical computers, as unitary 
transformations on qubits. It is known how in principle to do this, so 
that the quantum computer could compute anything which a classical 
computer can compute. The idea is to make use of the parallelism of 
complex superpositions, and entanglement between many qubits, to 
allow extremely fast algorithms for previously hard problems. Program 
and data, in the form of a sequence of binary digits, would be 
put into the quantum computer as the states $\ket 0$, $\ket 1$ of each of 
the component qubits. Then unitary evolution takes over in the 
product space, and leads after some time interval to a new joint 
state. The final step is to read out again, somehow, an output of 
the computation, and for that we must wait till the last ingredient 
has been discussed, measurement.

Already with just two qubits, entanglement can produce fascinating 
effects. In Section 3 we will use the entangled state of two qubits
$\frac 1 {\sqrt 2}\ket 0\otimes\ket 1 - \frac 1 {\sqrt 2}\ket 1\otimes\ket 0$ in order to 
perfectly teleport another quantum state from one location to another.

\subsection{Measurement}

So far we have described only the internal behaviour of quantum 
systems. Without any description of how the state of a quantum system can
have an influence on events in the classical outside world in which you 
and I walk about, and where we see tables and chairs, live or dead 
cars, not complex vectors or tensor products, the theory is completely empty. 
Moreover, so far the theory has been completely deterministic. Statisticians 
and probabilists will be getting impatient.

We describe here the most basic way in which we can obtain 
information from a quantum system. It is called a simple measurement. 
The idea is that we take the quantum system, bring it into interaction 
with some macroscopic experimental apparatus, and get to see some 
changes in the real world, which we quantify as a numerical 
measurement outcome $x$. The quantum system itself is changed by this 
process: one of the key ideas of quantum mechanics is that you cannot 
measure a system without disturbing it in some way. The process is 
random: both the outcome $x$ and the final state of the quantum 
system are random. But if for a given apparatus or experimental design,
we know the initial or input state $\psi$, and the 
outcome $x$, we also know the final or output state. The 
probability distribution of the outcome depends on the initial state, 
and on which of the many possible measurements---which of the possible 
experimental apparatusses---we use. 

The mathematical description goes as follows. Each measurement 
corresponds to a collection of orthogonal subspaces $A_{\{x\}}$ of our state 
space, labelled by the possible real values $x$ of the outcome. 
In our finite dimensional set-up there can be only a finite number of 
them, varying through some subset $\mathcal X$ of the real numbers.
The subspaces must not only be orthogonal but also span the whole state 
space, so that any state vector can be written as the sum of its 
orthogonal projections onto each of the subspaces $A_{\{x\}}$. Write 
$\Pi_{\{x\}}$ for the orthogonal projector onto $A_{\{x\}}$. Then applying 
the measurement described by $\{(x,A_{\{x\}}):x\in\mathcal X\}$
to a quantum system in state $\psi$ produces the 
value $x$ with probability $\|\Pi_{\{x\}}\psi\|^{2}$, the squared length 
of the projection of the state vector into the subspace $A_{\{x\}}$, and 
in this case the final state of the quantum system is just the 
renormalized projection $\Pi_{\{x\}}\psi/\|\Pi_{\{x\}}\psi\|$. By 
Pythagoras, and since we started with a normalized state vector, the 
sum of the squared lengths of the projections onto the orthogonal, 
spanning, subspaces $A_{\{x\}}$ equals $1$. And of course these squared 
lengths are real nonnegative numbers: thus, bona fide probabilities.
There is no harm in augmenting our collection of outcomes $\mathcal X$ 
with further values $x$ corresponding to $0$-dimensional subspaces 
$A_{\{x\}}$ 
consisting just of the zero vector. The length of the projection of 
$\psi$ onto the null subspace is zero, so this outcome is never 
observed. And a null subspace is orthogonal to every subspace.

An even more special case has each subspace $A_{x}$ one-dimensional, 
thus of the form $A_{\{x\}}=[\phi_{x}]$ for an orthonormal basis 
$\phi_{x}$ indexed by $x\in\mathcal X$. Then since $\Pi_{\{x\}}=
\ketbra{\phi_{x}}\phi_{x}$ one quickly sees that the result of the 
measurement is to yield the value $x$ with probability $|\braket 
{\phi_{x}}\psi|^{2}$, in which case the final state of the quantum 
system is $\ket{\phi_{x}}$. The complex numbers $\braket 
{\phi_{x}}\psi$ are called the probability amplitudes for the 
transition from $\psi$ to $\phi_{x}$, $x\in\mathcal X$.

There are a couple of alternative ways to mathematially reformulate 
this description. One way is to note that for a given measurement 
$\{(x,A_{\{x\}}):x\in\mathcal X\}$, $\mathcal X \subseteq\mathbb R$,
the collection of subspaces and values (except for null subspaces, 
which are irrelevant) can be recovered from the 
single matrix $X=\sum_{x\in \mathcal X} x \Pi_{\{x\}}$. This matrix 
is self-adjoint; it has real eigenvalues $x$ and its eigenspaces are 
the corresponding $A_{\{x\}}$. So the matrix $X$ is a compact mathematical 
packaging of all the information which we need to specify a measurement.
In physics such matrices are called observables, or physical 
quantities. Examples we have already seen are the Hamiltonian $H$, or 
for two-level systems, the spin observables 
$\sigma_{x}$, $\sigma_{y}$ and $\sigma_{z}$.

This compact mathematical formulation is moreover very powerful. 
Suppose we `measure the observable $X$' on the quantum system with 
state vector $\ket\psi$, state matrix $\rho=\ketbra\psi\psi$.
The probability to get the outcome $x$ is $\|\Pi_{\{x\}}\|^{2}
=(\Pi_{\{x\}}\psi)^{*}{\Pi_{\{x\}}\psi}
=\psi^{*}\Pi_{\{x\}}^{*}\Pi_{\{x\}}\psi=$
 (since a projection matrix is 
self-adjoint, $\Pi^{*}=\Pi$, and idempotent, $\Pi^{2}=\Pi$)
$=\psi^{*}\Pi_{\{x\}}\psi=$ (since a real number is a one-by-one 
matrix, hence equal to the trace of that matrix) 
$=\mathrm{trace}(\psi^{*}\Pi_{\{x\}}\psi)=$ (since one may cyclically 
permute matrix factors inside a trace) 
$=\mathrm{trace}(\psi\psi^{*}\Pi_{\{x\}})=\mathrm{trace}(\rho\Pi_{\{x\}})$.
Now multiply each probability by the value of the outcome $x$, and add 
over the values $x$; since $X=\sum_{x}x\Pi_{\{x\}}$ we find the 
celebrated trace rule, a most beautiful formula:
$\mathrm E_{\rho}(\mathrm{meas}(X))=\mathrm{trace}(\rho X)$
where $\mathrm{meas}(X)$ denotes the random outcome of measuring the 
observable $X$, and $\mathrm E_{\rho}$ denotes mathematical 
expectation when the (matrix) state of the quantum system is $\rho$.

This little formula: assigning a mean value under state $\rho$ to an 
observable $X$ (both represented by matrices, or in greater 
generality, operators), is the starting point for the field of 
quantum probability, which sees the mathematical structure of 
self-adjoint matrices (observables) and states, as analogous to the 
usual set-up of random variables and probability measures in 
classical probability theory. We have a way to compute expected values
$\mathrm{trace}(\rho X)$ somehow analogous to the classical formula 
(where now $X$ is a random variable on some probability space) 
$\int X\mathrm d\mathrm P$. I shall come back to this analogy, in the 
afterlude to the paper. However for us, the observable $X$ is just a 
convenient packaging of its eigenspaces and eigenvalues, and does not 
have an intrinsic role to play as a matrix or operator somehow acting 
on (multiplying) state vectors. But I would like to mention a further 
ramification. For a matrix $X=\sum_{x}x\Pi_{\{x\}}$ and a real 
function $f$, one can define the same function of the matrix $f(X)$ as
$f(X)=\sum_{x}f(x)\Pi_{\{x\}}$. Thus: keep the same eigenspaces, 
replace the eigenvalues by $f$ of the original eigenvalues. Well, 
this description is correct if the function $f$ is one-to-one; 
otherwise it should be modified to say: replace the eigenvalues by $f$ 
of the eigenvalues, if several eigenvalues map to the same function 
value, then merge the corresponding eigenspaces (i.e, take their 
linear span).  If the function $f$ is `square' or `exponential', then 
this curious definition does correspond to the existing, more orthodox 
definitions of $X^{2}$ or $\exp(X)$ for a given matrix $X$. Now given an
observable $X$ and a function $f$ we can talk about two different 
experiments: measure $X$ and evaluate the function $f$ on the outcome;
or directly measure $f(X)$. The resulting state of the quantum system 
is different if $f$ is many-to-one so that eigenspaces have merged;
one does not project so far when measuring $f(X)$ as with measuring $X$.
But it is a theorem that the probability distribution of the 
outcomes under the two scenarios is equal, and hence so are the 
expected values: $\mathrm E_{\rho}(f(\mathrm{meas}(X)))=
\mathrm E_{\rho}(\mathrm{meas}(f(X)))=\mathrm{trace}(\rho 
 f(X))$. I call this little formula, the law of 
the unconscious quantum physicist, since it is exactly analogous to 
the infamous law of the unconscious statistician in probability 
theory: the result that you can compute the expectation of a 
function $f$ of a random variable $X$ in two different ways: 
by integrating $f(x)$ with 
respect to the law of $X$, and by integrating $y$ with respect to the 
law of $Y=f(X)$. The quantum version of this law is part of the 
standard apparatus of quantum mechanics, and plays moreover a central 
role in foundational discussions, but is hardly ever explicitly stated 
let alone proved. Note that it leads to the computation not only of 
expectations but also of complete probability distributions: if I 
know the mean of every function of the outcome of measuring the 
observable $X$, I can recover the probability distribution of the 
outcome of measuring $X$. Thus all the expected values $\mathrm{trace}(\rho 
 f(X))$ do enable one to build a complete probability theory.

That was one way to mathematically reformulate measurement. Another 
way goes in the opposite direction: from relatively compact to 
over-elaborate. Yet it is very important for future developments 
(Section 4). For any measurable subset of the real line $B$, form the 
matrix $\Pi_{B}=\sum_{x\in B\cap\mathcal X}\Pi_{\{x\}}$. For each set 
$B$ this is a projection matrix, which projects into the sumspace of 
the eigenspaces of $X$, for $x\in B$. As such it satisfies the three 
axioms of a probability measure on the real line, but now with 
numbers replaced by matrices: (i), $\Pi_{B}\ge \mathbf 0$ for all $B$;
(ii), $\sum_{i}\Pi_{B_{i}}=\Pi_{B}$ for all disjoint $B_{i}$ with $B=\cup_{i}B_{i}$;
(iii), $\Pi_{\mathrm R}=\mathbf 1$. For a self-adjoint matrix $X$, the inequality 
$X\ge \mathbf 0$ means $\bra\psi X \ket \psi\ge 0$ for all vectors 
$\ket\psi$. We can now rewrite our probability rule for the 
probability distribution of the outcomes as 
$\mathrm P_{\rho}(\mathrm{meas}(X)\in B)=\mathrm{trace}(\rho \Pi_{B})$ 
for all measurable subsets $B$ of our outcome space (now considered to 
be all the real numbers). For the kind of measurements considered so 
far, the matrices $\Pi_{B}$ are not just self-adjoint and nonnegative, 
but also projection matrices: idempotent, as well. 
We call such a collection of matrices, a Projection-valued 
Probability Measure, or ProProM for short. In Section 4 we 
will see that it is necessary to take a wider view of measurement. 
There we will meet the notion of a generalized measurement, in which 
we replace the projection matrices $\Pi_{B}$ by arbitrary 
self-adjoint matrices, but still subject to the three rules of 
probability. We also call such generalized measurements, or rather 
their mathematical representations, Operator-valued Probabilty 
Measures or OpProM's.

In the previous section we formed the quantum analogue of product 
probability spaces. Also observables and measurements on one 
component of a product system can be considered as defined on a 
product system. The observable $X$ on a subsystem corresponds to the 
observable $X\otimes\mathbf 1'$ on the product of that system with 
another: same eigenvalues, eigenspaces equal to the original 
eigenspaces tensor product with the other complete space. If $X$ and 
$Y$ are two observables of two different quantum systems then 
$X\otimes 1'$ and $\mathbf 1 \otimes Y$ give an example of commuting 
observables: their product, taken in either order, is the same. 
Observables which commute model measurements which may be done 
simultaneously. Whether one first measures the one, then the other, 
or vica-versa, the probabilistic description of joint outcome and of 
final state is identical. A product system is often used to model a 
pair of particles at two different locations, and the observables of 
each subsystem correspond to measurements which may be made at the 
two separate locations, and which naturally do not influence one 
another's outcome. In particular, if the product system is in a 
product state, then the outcomes of measurements on the subsystems are
independent with the same distribution as if everything had been 
considered separately, as one naturally would desire.

\subsubsection{Example of measurement: the qubit}

For a $2$-dimensional state space 
one can only find sets of pairs of non-trivial, orthogonal subspaces, 
each pair corresponding to a pair of othonormal basis vectors. Now as 
we sketched previously, there is a one-to-one correspondence between 
state-vectors of $\mathrm C^{2}$ and directions (unit vectors) in 
$\mathrm R^{3}$. Orthogonal state vectors correspond to opposite 
directions. Let us label the two possible outcomes of one of these 
measurements, by the real values $+1$ and $-1$. Then each of the 
non-trival simple measurements corresponds to the two projector 
matrices $\ketbra{\vec v}{\vec v}$ and $\ketbra{-\vec v}{-\vec v}$.
The two add up to the identity matrix, and can also be written, as we 
saw before as $\frac12(\mathbf 1\pm\vec v\cdot{\boldsymbol \sigma})$.
The corresponding observable (matrix) is 
$\ketbra{\vec v}{\vec v}-\ketbra{-\vec v}{-\vec v}= \vec v\cdot{\boldsymbol 
\sigma}$,
or as the physicists say `the spin observable in the direction $\vec v$'.
A little computation shows that the probability of the two outcomes 
$\pm1$,
when this observable is measured on a quantum system in the state 
$\ket{\vec u}$, is $\frac12(1\pm\vec u\cdot\vec v)$. The resulting 
state of the particle is $\ket{\pm\vec v}$. This has the implication 
that one can prepare particles in a given state, say $\ket{\vec 
v}$,by measuring particles in any state and only keeping those, for 
which the outcome was $+1$. Thus measurement, often thought of as 
being a final stage of an experiment, might also be the initial stage 
called `preparation'.

This measurement is realized on the spin of electrons in a so-called 
Stern-Gerlach device, a specially shaped magnet which can be 
physically oriented in the real direction $\vec v$ and carries out 
precisely the measurement just described. Electrons leave the magnet 
in two streams, in one stream all particles have the state  $\ket{\vec 
v}$, in the other they all have the state $\ket{-\vec v}$. The 
relative sizes of the two output streams depends on the initial states 
of the electrons.

\section{INTERMISSION: THE EXAMPLE OF TELEPORTATION}

We will illustrate the ingredients by the beautiful example of 
quantum teleportation, discovered by Charles Bennett (IBM) {\it et 
al.}\ in the mid nineties, and done in the laboratory, just a couple of years later, 
by Anton Zeilinger, in Innsbruck. Since then the experiment has been repeated 
in many places. The experiment is done with polarized photons, and the 
basis states can be thought of as $\ket\leftrightarrow$ 
($x$ direction), $\ket\updownarrow$ ($y$-direction). 

It is useful here to give some further discussion of how polarization 
of photons can be reformulated in the language of qubits.
Think of light coming 
towards you in the $z$ direction, and oscillating sinusoidally, 
with the same frequency, but possibly different relative amplitude and phase,
in both both the $x$ direction and the $y$ direction. The oscillations 
generate a (perhaps flattened) spiral around the $z$ direction, coming 
towards you.
Head on, you see an elliptical motion around the $z$ axis which might 
be directed clockwise or anticlockwise; the ellipse might be perfectly 
circular or 
perfectly flat (a line segement) or anything in between; the 
orientation of the major axis of the ellipse can be anything in the 
$x$-$y$ plane.
The perfectly flat version is how light comes 
out of a polarization filter (e.g., your sunglasses: the oscillation 
occurs entirely in one plane). Now imagine mapping all the different 
`directed, oriented, ellipses' onto the surface of the three-dimensional real sphere as 
follows: the clockwise ellipses on the Northern hemisphere, the 
anticlockwise on the Southern; the `flat' ellipses are arranged 
around the equator, and the two circles are at the North Pole and the 
South Pole. As one moves completely around the earth, at constant latitude, the 
direction of the ellipse rotates slowly around $180^{\circ}$. In short: 
all possible polarizations of light (all possible shapes of directed, 
oriented, ellipses) 
can be mapped one-to-one onto the directions in real three-dimensional space.

Now light behaves both as a wave and as a stream of particles 
(photons). In fact this is the essence of a quantum mechanical 
description; what we now know is that wave-particle duality extends to 
all known physical objects (for instance: photons, electrons, neutrons, protons; 
but also at higher and lower scales).
The quantum state of polarization of one photon is described by 
a two dimensional state vector $\ket{\vec u}$.  
All possible transformations of the state of polarization
correspond to orthogonal rotations of the real vector $\vec u$, and 
to unitary transformations of the quantum state vector $\ket{\vec u}$.
They can be implemented in the laboratory by passing the light 
through suitable transparent media (fluids and crystals). Moreover 
any simple measurement or preparation can be implemented with beam splitters and 
polarization filters.

Now the problem of teleportation is as 
follows. Alice, who lives in Amsterdam, 
is given a qubit (polarized photon) in an unknown 
state, say $\alpha\ket\leftrightarrow+\beta\ket\updownarrow$. She wants 
to transmit it to Bob, who lives in Beijing, 
and she can only communicate with Bob by email. 
(If you prefer, replace Amsterdam and Beijing with, perhaps 
futuristically,  P'yongyang and Seoul). What can 
she do? She could measure the qubit, e.g., look to see if the photon is polarized 
$\leftrightarrow$ or $\updownarrow$. She gets the answer: ``$\leftrightarrow$'' or 
``$\updownarrow$''; the answer is random, with 
probabilities $|\alpha|^{2}$, $|\beta|^{2}$ depending on the unknown $\alpha$, $\beta$. 
The photon's original state is 
destroyed, we cannot learn anything more about it. So all she could do is email to 
Bob: ``I saw (e.g.) $\leftrightarrow$''. He makes a horizontally polarized photon. This is a 
poor, random, copy of the original one, and the original one has gone. Can they 
do better? Well, there are many other measurements Alice could make, but they 
all have the same property, of only providing a small, random, amount of 
information about the original state, and destroying it in the process. In fact it is 
a result from the theory of quantum statistical inference due to
\citet{helstrom67, helstrom76, braunsteincaves94}, that whatever measurement is 
carried out by Alice, the Fisher 
information matrix based on the probability distribution of the outcome of the 
experiment, concerning the unknown parameters $\alpha,\beta$,  has a strictly positive 
lower bound. The famous no-cloning theorem could also be invoked here: 
it is impossible to convert one quantum system into two identical 
copies. We will review this result in Section 5.

In order to succeed, Alice and Bob need a further resource. What they do is 
arrange that each of them has another photon, these two (extra) photons in the 
entangled joint state  $\frac 1{\sqrt 2}\ket 0 \otimes \ket 1   - 
\frac 1{\sqrt 2} \ket 1  \otimes 
\ket 0$. This particular state is called the singlet, or Bell state.
This is nowadays a routine matter. It is created by having someone 
else, at a location between Amsterdam and Beijing, excite a Calcium 
atom with a laser in such a way that the atom moves to a higher energy 
level. Then the energy rapidly decays and two photons are emitted, in 
equal and opposite directions. One travels to Amsterdam, the other to 
Beijing. Now we have three qubits, living together in 
an eight-dimensional space, of which four of the dimensions---two of the 
qubits---are on Alice's desk, the other two dimensions---one qubit---on 
Bob's desk. Below we will see three lines of elementary algebra, with the astounding 
implication that Alice can carry out a measurement on her desk, get one of $4$ 
random outcomes, each with probability $1/4$, then email to Bob which outcome 
she obtained; he correspondingly carries out one of $4$ different, prescribed, unitary 
operations, and now his photon is magically transformed into an identical copy 
of the original, unknown, qubit which was given to Alice. Two (unknown) 
complex numbers $\alpha$ and $\beta$ have been transmitted, with complete accuracy, by 
transmitting two bits of classical information. (More accurately, two 
real numbers, say $(\theta,\phi)$; but this is just as amazing).

Now it is worth asking: how can we know that a certain experiment has actually 
succeeded? The answer is of course by statistics. One needs, many times, to 
provide Alice with qubits in various states. Some of these times, the qubits are not 
teleported, but are measured in Alice's laboratory. On the other occasions, the 
qubits are teleported to Bob, and then measured in Bob's laboratory. The 
predictions of quantum theory are that the {\it statistics\/} of the measurements at 
Alice's place, are the same as the statistics of the measurements at Bob's place.

So suppose a single spin-half particle with state-vector 
$\alpha\ket 0+\beta\ket 1$  is brought into interaction with a 
pair of particles in the singlet state, written abbreviatedly as 
$\ket{01}-\ket{10}$ (and discarding a constant factor). I am using the following 
shorthand: for instance, $\ket 0 \otimes \ket 0 \otimes \ket 1 $ is 
written as $\ket {001}$. The order of the three components 
is throughout: first the particle to be teleported (on Alice's desk 
in Amsterdam), then 
Alice's part of the singlet pair (also on her desk), then Bob's part 
of the singlet pair (on his desk in Beijing). The whole $2^{3}$ 
dimensional system has state-vector, 
multiplying out all (tensor) products of sums of state vectors, 
and up to a factor $1/\sqrt 2$,
$(\alpha\ket 0+\beta\ket 1)\otimes(\ket{01}-\ket{10})=
\alpha\ket{001}-\alpha\ket{010}+\beta\ket{101}-\beta\ket{110}$.
Now we introduce the following four orthogonal state-vectors for the two 
particles in Amsterdam, neglecting another constant factor $1/\sqrt 2$, 
$\Phi_1=\ket{00}+\ket{11}$, 
$\Phi_2=\ket{00}-\ket{11}$,
$\Psi_1=\ket{01}+\ket{10}$, 
$\Psi_2=\ket{01}-\ket{10}$,
and we note that our three particles together are in a pure state 
with state-vector which may be
written (up to yet another factor, $1/\sqrt 4$)
$\alpha(\Phi_{1}+\Phi_{2})\otimes\ket 1
-\alpha(\Psi_{1}+\Psi_{2})\otimes\ket 0
+\beta(\Psi_{1}-\Psi_{2})\otimes\ket 1
-\beta(\Phi_{1}-\Phi_{2})\otimes\ket 1
$.
Rearranging these terms (noting that $\alpha$ and $\beta$ are numbers 
so can be moved through the tensor products at will) one finds the 
state
$\Phi_1\otimes(\alpha\ket 1-\beta\ket 0)
+\Phi_2\otimes(\alpha\ket 1+\beta\ket 0)
+\Psi_1\otimes(-\alpha\ket 0+\beta\ket 1)
+\Psi_2\otimes(-\alpha\ket 0-\beta\ket 1)$.
So far nothing has happened at all: we have simply rewritten the 
state-vector of the three particles as a superposition of four 
state-vectors, each lying in one of four orthogonal two-dimensional
subspaces of $\mathbb C ^2 \otimes \mathbb C ^2 \otimes \mathbb C ^2$:
namely the subspaces
$[\Phi_{1}]\otimes\mathbb C ^2$,
$[\Phi_{2}]\otimes\mathbb C ^2$,
$[\Psi_{1}]\otimes\mathbb C ^2$ and
$[\Psi_{2}]\otimes\mathbb C ^2$.

To these four subspaces corresponds a simple measurement. 
It only involves the two particles in Amsterdam and hence may be carried out by
Alice. She obtains one of four different outcomes, each with probability 
$\frac14$, so she learns nothing about the particle to be teleported.
However, conditional on the outcome of her measurement, 
the particle in Beijing is in one 
of the four states 
$\alpha\ket 1-\beta\ket 0$,
$\alpha\ket 1 +\beta\ket 0$, 
$-\alpha\ket 0 +\beta\ket 1$,
$-\alpha\ket 0-\beta\ket 1$. 
So Bob knows that he has with probability $\frac14$, either of those
four states. It can be verified that whatever he does with that 
particle, his statistical predictions are the same as before Alice 
made her measurement: nothing has changed at Beijing, yet!
But once the outcome of the measurement at Amsterdam is transmitted to 
Beijing (two bits of information, transmitted by classical means),
Bob is able by means of one of four unitary transformations to 
transform the resulting pure state into the state with state-vector 
$\alpha|0\rangle+\beta|1\rangle$. For instance, if the first of the 
four possibilities is realized, Bob must change the sign and carry out 
a spin-flip to convert $\alpha\ket 1-\beta\ket 0$ into $\alpha\ket 0+\beta\ket 1$.
He does not need any knowledge of $\alpha$ and $\beta$ to do this: he 
just carries out two fixed unitary transformations. In each of the 
four cases, there is a fixed unitary transformation which does the job.

Neither Alice nor Bob learn anything 
at all about the particle being teleported by this procedure. In fact, 
if they did get any information about $\alpha$ and $\beta$ the 
teleportation would have been less than succesful. One cannot learn 
about the state of a quantum system without (at least) partially 
destroying it. The information one gains is random. There is no going 
back.

\section{MODEL GENERALIZATION AND SYNTHESIS}

We are close to describing new and interesting statistical problems. 
However, first we must extend the notion of state, and the notion of 
measurement, used so far. Suppose we want to get information about the 
state of some quantum system. There is more that we can do than just 
carry out one simple measurement on the given system. We could for 
instance first bring the system being studied into interaction with 
another quantum system, in some known state. After a unitary 
evolution of the joint system, one could measure the auxiliary system.
Next, discard this system, and bring the original particle (which is 
now in some new state, dependent on the results so far), into 
interaction with another auxiliary system. Do the same again. At each 
stage one could allow the various operations (initial state 
of the auxiliary system, unitary transformation, measurement of 
auxiliary system \ldots) to depend on the outcomes obtained so far. 
Finally after some number of operations, take some function of all 
the outcomes obtained in the intermediate steps. 

This provides a vast repertoire of possible strategies, and it seems impossible to 
describe ``everything that can be done'' in a concise and 
mathematically tractable way, in order to optimize over this collection.
It is not actually clear in advance that these more elaborate 
measurement schemes could be useful, but it is a fact that they arise 
in practice, and moreover that they often provide strictly better 
solutions to statistical design problems, than the simple 
measurements!

Secondly, the notion of state is just a little restricted. Suppose 
that each time a qubit is manufactured, slight variations in 
temperature, materials, and so on produce slightly different states. 
The identical copies we are given are not single qubits in an 
elementary, so-called pure, 
state, but are actually i.i.d.\ drawings from a probability distribution over 
pure states. It seems that we need to know: the complete distribution of 
pure states, that the experimenter is sampling from. Again this would seem to be 
an unwieldly, complicated, object. 

Amazingly both the complications lead after a beautiful synthesis 
into generalized notions of state and of measurement which are very compact 
and amenable to mathematical analysis. Moreover, the syntheses 
(very different composite measurements may be represented by the 
same, compact, mathematical object, and similarly, completely different 
probability distributions over pure states cannot be distinguished 
either) highlights new and extraordinary 
features of quantum reality.

Recall that we could describe the state of a quantum system with the 
matrix $\rho=\ketbra\psi\psi$ rather than the vector $\ket\psi$. Now 
suppose that according to one scenario,
quantum systems in states $\ket{\psi_{i}}$ are produced with 
probabilities $p_{i}$, and measured in any complicated way allowed by 
the rules of quantum mechanics (i.e., using the ingredients of 
Section 2, in any combination). In another scenario, quantum systems 
in states $\ket{\phi_{j}}$ are produced with 
probabilities $q_{j}$, and measured in the same way. 
Suppose that the two scenarios are such the \emph{average} state matrix is 
the same: 
$\sum_{i}\ketbra{\psi_{i}}{\psi_{i}}=\sum_{j}\ketbra{\phi_{j}}{\phi_{j}}=\rho$, say. 
Suppose the final outcome of measurement is some outcome $x$ in an 
arbitrary (now possibly very large) measurable sample space $(\mathcal 
X,\mathcal B)$. Now if 
we have specified the measurement procedure, however complicated it 
is, then the rules from Section 2 allow us to compute the probability 
law of the random outcome $X$ under either of the two scenarios. Then 
one can state the following theorems:

\begin{thm}The probability distribution of $X$, i.e., the collection 
of probabilities
    $(\Pr(X\in B):B\in\mathcal B)$, only depends on the 
average or mixed state $\rho$; i.e., it is the same under our two 
scenarios, whatever the measurement protocol. Moreover the 
mapping from mixed state matrix $\rho$ to probability law of outcome 
is affine, i.e., linear under convex combinations of state matrices 
$\rho$.
\end{thm}

\begin{thm}Any affine mapping from mixed state matrices $\rho$ to 
probability distributions on $(\mathcal X,\mathcal B)$
is of the form $\mathrm P_{\rho}(X \in 
B)=\mathrm{trace}(\rho M(B))$ where 
$(M(B):B\in\mathcal B)$ is an Operator-valued Probability Measure 
(OProm); i.e., $M(B)$ is a self-adjoint matrix for every $B$ 
satisfying the axioms of a probability measure: $M(B)\ge \mathbf 0$ for all 
$B$;
$B(\mathcal X)=\mathbf 1$; $M(B)=\sum_{i}M(B_{i})$ whenever $B$ is 
the disjoint countable union of $B_{i}$.
\end{thm}

\begin{thm}Any operator-valued probability measure can be realized by 
bringing the quantum system being measured into interaction with an 
auxiliary system (so-called ancilla) in some fixed state $\rho_{0}$, 
applying a unitary evolution to the joint system, applying a simple 
measurement to the ancilla, and discarding the ancilla.
\end{thm}

This sequence of results tells us: everything that is allowed by 
quantum mechanics, is necessarily of the form of an OProM. And 
conversely, every OProM can in principle be realized, by a procedure 
which one might call quantum randomization, since it is based on 
forming a product system with a completely independent system, and 
then measuring the joint system. (In the literature, the abbreviation 
POM or POVM is often used, standing for `probability-operator 
measure', or `positive operator valued measure'; in our opinion that 
nomenclature is inaccurate).

Every mixture of state matrices is a non-negative self-adjoint matrix 
with trace $1$. Such a matrix is called a density matrix and every 
density matrix can be realized as a mixture of pure states, states of 
the form $\ketbra\psi\psi$, in general 
in very many ways. The pure states have density matrices which are 
idempotent, $\rho^{2}=\rho$. These cannot be written as a probability 
mixture over more than one state. Recall that the simple measurements 
could be represented with Projection-valued Probability Measures. So a 
modest mathematical extension of our basic notions allows us to 
encompass everything that quantum mechanics allows, in a concise and 
powerful way.

The underlying 
mathematical theorems here are due to Naimark, Holevo, Ozawa and others;
see \citet{helstrom76}, \citet{holevo82}.
They can be extended to describe in precisely the same way, not just the 
mapping from input state to observed data, but also to observed data 
and output state, conditional on the observed data.  This leads 
to the somewhat sophisticated mathematical notion of completely positive 
instruments and conditional states; the main theorems are due to
Stinespring, Davies, Kraus and again Ozawa. The paper 
\citet{barndorffnielsenetal01} contains many references to these and further 
devlopments. In particular there is great interest 
presently in modelling continuous time measurement of a quantum 
system, or continuous time interaction of a quantum system with a 
much larger environment, leading to a rich theory of quantum 
stochastic processes.

\subsubsection{Example: the qubit}

Recall that the pure state matrices of a qubit are of the form 
$\frac12(\mathbf 1 +\vec u\cdot\vec{\boldsymbol\sigma})$ where $\vec u$ 
is a unit vector in real three-space. Arbitrary probability mixtures 
of such states (corresponding to preparing a pure state 
chosen by a classical randomization from some probability distribution 
over the unit vectors in $\mathrm R^{3}$) can therefore be completely 
described by the resulting mixture of state matrices, which must be of the 
form $\rho=\frac12(\mathbf 1 +\vec a\cdot\vec{\boldsymbol\sigma})$ where 
now $\vec a$ is an arbitrary vector in the real three-dimensional unit 
ball. A simple measurement of spin in the direction $\vec v$, of this 
quantum system, results in outcomes $\pm 1$ with probabilities
$\frac12(1\pm\vec v \cdot \vec a)$.
If we had many copies of the quantum system, we could determine the 
vector $\vec a$ to arbitrary precision by carrying out large numbers 
of measurements of spin in three orthorgonal directions, e.g, the 
$x$, $y$ and $z$ directions. Is this the most accurate way to 
determing $\vec a$ when we have a large number $N$ of copies at our 
disposal?

The generalized measurements or OProM's form a huge class of possible 
experimental designs. Here we just mention one such measurement.
It has an outcome space consisting of just three outcomes, let us 
call them $1$, $2$ and $3$. Let $\vec v_{i}$, $i=1,2,3$, denote three 
unit vectors in the same plane through the origin in $\mathrm R^{3}$, 
at angles of $120^{\circ}$ to one another. Then the matrices 
$M(\{i\})=\frac 13 (\mathbf 1 +\vec 
v_{i}\cdot\vec{\boldsymbol\sigma})$  define an operator-valued 
probability measure on the sample space $\{1,2,3\}$ which is called 
the triad, or Mercedes-Benz. It turns up as the optimal solution to 
the decision problem: suppose a qubit is generated in one of the 
three states $\ket{\vec v_{i}}$, $i=1,2,3$, with equal probabilities. 
What decision rule gives you the maximum probability of guessing the 
actual state correctly? There is no way to equal the success 
probability of this method, if one only uses simple measurements, even 
allowing for (classically) randomized procedures. One could say that 
quantum randomization is sometimes necessary to maximally extract information 
from a quantum system. The triad could be realized by bringing the 
system under study into interaction with another three-dimensional 
system in a certain, fixed, state, carrying out a certain unitary 
transformation on the joint system, and then carrying out a certain 
simple measurement on the ancilla.

Another measurement which occurs as the optimal solution to some 
estimation problems has outcomes which are continuously distributed 
real unit vectors; the matrix elements of the OProM $M(B)$ have density 
$\frac1{4\pi}(\mathbf 1+\vec v\cdot \vec{\boldsymbol\sigma})$ with 
respect to Lebesgue (surface) measure on the unit sphere. It would be 
realized in practice by coupling the qubit to a quantum system with 
infinite dimensional state space.

\section{QUANTUM STATISTICS: DESIGN AND INFERENCE}

Suppose we are given $N$ qubits in an identical, unknown, state, what is the 
best way to determine that state? It is known (by the  statistical information 
bound we are about to discuss) that whatever one does, one cannot achieve better than a 
certain degree of accuracy, of the order of size of $1/\sqrt N$. It is not known 
what constant over $\sqrt N$, is best. And a most intriguing question, only partially 
solved, is: does it pay off to consider the $N$ qubits as one joint system, having a 
state of a the special form $\rho^{(N)}=\rho^{\otimes N}$ in a $2^N$ dimensional state space, or can one just as well 
measure them separately? Note that by considering the $N$ copies as 
one collective system, we have a much 
vaster repertoire of possible measurements, so from a mathematical point of 
view, the answer should surely be that joint measurements pay off. However 
physical intuition would perhaps say the opposite. I have worked on asymptotic 
versions of this problem. So far physicists have hardly considered this route, and 
the literature has mainly seen calculations in rather special situations ($N=2$, for 
instance), with conclusions which depend on all kinds of features of the 
problem---prior distributions, loss functions---which are really arbitrary. The advantage of 
my approach is that these extraneous and arbitrary features become irrelevant for 
large, but finite $N$; the problem {\it localizes}, second order approximations are good, 
loss functions might as well be quadratic, prior distributions are irrelevant. Using 
the van Trees inequality (a Bayesian Cram\'er--Rao bound, see 
\citet{gilllevit95}) I have, together with 
Serge Massar, derived {\it frequentist\/} large sample results on what is asymptotically 
best, under various measurement scenarios; see the survey paper
\citet{gill01} and the original work 
\citet{gillmassar00}. Further results are contained in 
\citet{barndorffnielsenetal01}; and a more comprehensive survey paper 
by \citet{barndorffnielsenetal01b} is in preparation.

Similar results have been obtained, interestingly, with quite different 
methods, in a series of papers, by \citet{young75}, 
\citet{fujiwaranagaoka95}, \citet{hayashi97}, and 
citet{hayashimatsumoto98}.

The most exciting result we have found is as follows: if the unknown state is 
known to be pure, then a certain very simple but adaptive strategy of basic yes/no 
measurements on the separate qubits, achieves the maximal achievable accuracy. If 
however the state is mixed, then we do not know the best strategy. Limited to 
separate measurements, we do know what can be achieved. We know that joint 
measurements can achieve startling increases in accuracy. But we do not know 
how much can be maximally achieved (there are known bounds, but they are 
known to be unachievable). This seems to be a promising research direction.

The `pure state' solution is as follows. First get a rough estimate of 
the direction of spin by measuring the spin in the $x$, $y$ and $z$ 
directions separately, on a large number, but small fraction, of the 
particles; say, on $\sqrt N$ particles each. Now do a simple 
measurement of the spin on each half of the remaining $N-3\sqrt N$ 
particles, in two perpendicular directions \emph{orthogonal} to the direction of 
the rough estimate. In the physics literature it has been suggested 
that one should try as well as possible, to measure in the same 
direction as the unknown spin---basically the opposite to our solution. 
And the simple strategy just described, is asymptotically as good as 
anything else one can imagine, however sophisticated, on all $N$ 
particles together. In particular it is asymptotically as good as the 
the theoretically optimal solution for a uniform 
prior distribution and certain rather special loss functions, namely 
a beautiful but practically impossible to implement generalized
measurement on the collective of particles.

\subsection{Finite sample optimal design: the quantum information bound}

In this subsection I want to prove and discuss a central and now 
classical result on the design of optimal quantum measurements, the quantum Cram\'er--Rao 
inequality and quantum information bound. The quantum information 
matrix plays a key role in the results I have just mentioned, though 
new quantum information bounds are needed, as we will see.

We first introduce analogues to the score function and information 
matrix of classical statistics: the quantum score and the quantum 
information. Just as the classical score function can be thought of 
both as a random 
variable, and as the derivative of the logarithm of the probability 
density, so is the quantum score both an observable (self-adjoint 
matrix) and a certain kind of derivative of the density matrix. The
quantum information is the mean of the squared quantum score, just as 
in classical statistics, except that now the mean is taken using the 
trace rule for expectations of outcomes of measurements of observables.

Consider a quantum statistical model: that is 
to say a parametric family of density matrices 
$(\rho(\theta):\theta\in\Theta)$. 
A measurement $M$ with outcome space $(\mathcal X,\mathcal B)$
and with density $m$ with respect to a (real) sigma-finite 
measure $\mu$ is given. When we apply the measurement to a 
quantum system with state $\rho(\theta)$ in this model, we obtain an outcome
with density 
$p(x;\theta)=\mathrm{trace}(\rho(\theta)m(x))$ with respect to $\mu$.
For this classical parametric statistical model, one can compute the 
Fisher information matrix; we denote it as $I(\theta;M)$.

For the moment, suppose that the parameter space is one-dimensional.
We define the so-called quantum score as follows: it is implicitly defined as 
the self-adjoint matrix $\lambda=\lambda(\theta)$ which solves the 
equation $\rho'=\frac12(\lambda\rho+\rho\lambda)$. Here, $\rho'$ 
denotes the derivative of $\rho(\theta)$ with respect to $\theta$ (the 
matrix of derivatives of matrix elements). Just as the state $\rho$ 
depends on $\theta$, so also do $\rho'$ and $\lambda$. Now the quantum 
information (number) is defined as 
$I_{Q}(\theta)=\mathrm{trace}(\rho(\theta)\lambda(\theta)^{2})$. From 
what we learnt before, this number is the mean value of the square of 
the outome of a measurement of the observable $\lambda(\theta)$.
If the parameter $\theta$ is actually a vector, then one defines 
quantum scores component-wise, and finally defines the quantum 
information matrix elementwise by 
$I_{Q}(\theta)_{ij}=\mathrm{trace}
(\frac12\rho(\theta)(\lambda(\theta)_{i}\lambda(\theta)_{j}
+\lambda(\theta)_{j}\lambda(\theta)_{i}))$.

The following quantum information inequality due to 
\citet{braunsteincaves94} is crucial:
\begin{equation*}
I(\theta;M)\le I_{Q}(\theta)
\end{equation*}
for all measurements $M$. From this inequality one immediately has 
the quantum Cram\'er--Rao inequality, \citet{helstrom67}: for all measurements $M$, and 
any unbiased estimator $\widehat\theta$ based on the outcome of that 
measurement,
\begin{equation*}
\mathrm{Var}(\widehat\theta)\ge I_{Q}(\theta)^{-1}.
\end{equation*}

To prove the information inequality we need to express the Fisher information in 
the outcome of $M$ in terms of the quantum score. Since 
$p(x;\theta)=\mathrm{trace}(\rho(\theta)m(x)$ it follows that 
$p'(x;\theta)=\mathrm{trace}(\rho'(\theta)m(x))=\frac12(\mathrm{trace}(\rho\lambda m)
+\mathrm{trace}(\lambda\rho m))=\frac12(\overline{\mathrm{trace}((\rho\lambda 
m)^{*})}
+\mathrm{trace}(\lambda\rho m))=\frac12(\overline{\mathrm{trace}(m\lambda \rho)}
+\mathrm{trace}(\rho m\lambda))=\frac12(\overline{\mathrm{trace}(\rho 
m\lambda m)}
+\mathrm{trace}(\rho m \lambda ))=\Re(\mathrm{trace}(\rho m \lambda))$.
Thus the classical score function is $\Re(\mathrm{trace}(\rho(\theta) 
m(x) \lambda(\theta)))/p(x;\theta)$.

From now, $\theta$ is fixed. Define $\mathcal X_{+}=\{x:p(x;\theta)>0\}$ and $\mathcal 
X_{0}=\{x:p(x;\theta)=0\}$. 
Define $A=A(x)=m(x)^{\frac12}\lambda\rho^{\frac12}$,
$B=B(x)=m(x)^{\frac12}\rho^{\frac12}$, and $z=\mathrm{trace}\{ A^* B\}$.
Note that $p(x;\theta)=\mathrm{trace}\{ B^* B\}$.

The proof of the quantum information inequality 
consists of three steps.
The first will be an application of the trivial inequality
$\Re(z)^2\le |z|^2$ with equality if and only if $\Im(z)=0$. 
The second will be an application of the Cauchy--Schwarz inequality 
$| \mathrm{trace} \{A^* B\} |^2 \le \mathrm{trace} \{A^* A\} 
\mathrm{trace} \{B^* B\} $ 
with equality if and only if
$A$ and $B$ are linearly dependent over the complex numbers.
The last step consists of replacing an integral of a nonnegative function
over $\mathcal X_{+}$ by an integral over $\mathcal X$. Here is the 
complete proof:
\begin{eqnarray}
  I(\theta;M)  ~ &=& ~  \int_{\mathcal X_{+}}
        p(x;\theta)^{-1} ( \Re\,\mathrm{trace} (\rho\lambda m(x) )^2 \mu(\mathrm d  x) 
               \nonumber \\
  ~ &\le& ~ \int_{\mathcal X_{+}} p(x;\theta)^{-1} |\mathrm{trace} 
  (\rho\lambda 
               m(x))|^2  \mu(\mathrm d  x) 
                                                                \nonumber \\
  ~ &=& ~  \int_{\mathcal X_{+}} \left| 
         \mathrm{trace}  \left( \,  m(x)^\frac12  \rho^\frac12  )^* \,
                  ( m(x)^\frac12  \lambda \rho^\frac12 \, \right)
   \right|^2 (\mathrm{trace} (\rho m(x)))^{-1}\mu(\mathrm d  x)   \nonumber  \\
  ~ &\le& ~  \int_{\mathcal X_{+}} 
  \mathrm{trace} ( m(x) \lambda \rho \lambda ) 
                \mu(\mathrm d  x)  \nonumber \\
  ~ &\le& ~ \int_{\mathcal X}\mathrm{trace} ( m(x) \lambda 
             \rho \lambda ) \mu(\mathrm d  x)  \nonumber \\
  ~ &=&  ~  I_{Q}(\theta) .\nonumber
\end{eqnarray}
In the last step we used that $\int m(x)\mu(\mathrm d x)=M(\mathcal 
X)=\mathbf 1$.
One can verify that equality holds, if and only if 
$m(x)^{\frac12}\lambda(\theta)\rho^{\frac12}(\theta)
=r(x,\theta)m(x)^{\frac12}\rho^{\frac12}(\theta)$ for some real 
$r(x;\theta)$, for
$p(x;\theta)\mu(\mathrm d x)$ almost all $x$. Under smoothness, 
positivity and nondegeneracy conditions, this tells us that for optimal 
Fisher information, an attaining measurement $M$ can be nothing else than 
the simple measurement of the quantum score observable 
$\lambda(\theta)$. In general, this measurement does attain the 
information inequality.

For vector parameters the information inequality and Cram\'er--Rao 
inequality remain true; the proof follows by considering all smooth
one-dimensional submodels, whose classical and quantum score functions
are of course linear combinations of the component scores.

For the models for one 
qubit which we have been studying, the parameter $\theta$ might be 
taken to be the real vector $\vec u$ or $\vec a$ of a completely 
unknown pure state, or completely unknown mixed state. The quantum 
information matrices for either of these models is easy to compute, 
and is as one might expect strictly positive.

Now this result already tells us a great deal. First of all, it is 
not difficult to show that the quantum information for $N$ identical copies of a 
quantum statistical model, i.e., with density matrix 
$\rho^{(N)}(\theta)=\rho(\theta)^{\otimes N}$,  is $N$ times the 
information in one copy. Thus even if one uses elaborate measurements 
on a joint system of $N$ identical copies, one cannot beat the $\frac 1 {\sqrt N}$
of classical statistics. Moreover, just thinking about one copy: one 
cannot determine the quantum state exactly by doing elaborate enough 
measurements: otherwise the quantum information would not be strictly 
positive. And we have a proof of the no-cloning theorem: if by 
combining the basic ingredients of quantum mechanics in some way we 
could convert one copy of an unknown quantum state into two identical 
copies, we could make an arbitrary large number of identical copies, 
and hence estimate the state arbitrarily well, but this contradicts 
the positive information bound.

Much more comes out of it. Suppose the parameter is one-dimensional. 
Then the best measurement in terms of Fisher information is to 
measure the quantum score. But that typically depends on $\theta$ and 
moreover for different $\theta$, the scores $\lambda(\theta)$ do not 
commute. So there typically is no single measurement which achieves 
the information bound uniformly in the parameter value. However for 
large $N$ one can get close: using a small number of copies, get a 
rough estimate of $\theta$, then measure the `estimated score' on the 
remaining copies. For large $N$ this will be close to measuring the 
true but unknown score on all copies, hence close to attaining the 
information bound on the collective. And thus the maximum likelihood 
estimator based on the data, will approximately achieve the 
Cram\'er-Rao bound.

But now suppose the parameter is not scalar. Typically the  
score observables for the different components of $\theta$ do not commute.
This means that even if you (roughly) know the value of $\theta$, 
completely different and mututally incompatible experiments are 
needed to determine the different components of $\theta$ as well as 
possible. Now we seem to be stuck. The quantum information inequality 
is the best matrix inequality one can have, but it does not delineate 
the class of attainable classical information matrices; i.e., not 
every matrix $J$ with $J\le I_{Q}(\theta)$ is the information matrix 
of some measurement $M$ at $\theta$. Thus we need something better, in 
order to describe what can be done.

\subsection{Quantum Asymptotics}

In classical statistics, the Cram\'er--Rao bound is attainable 
uniformly in the unknown parameter only under rather special circumstances.
On the other hand, the restriction to unbiased estimators is hardly 
made in practice and indeed is difficult to defend. 
However, we have a richly developed asymptotic theory 
which states that in large samples certain estimators (e.g., the 
maximum likelihood estimator) are approximately unbiased and 
approximately normally distributed with variance attaining the 
Cram\'er--Rao bound. Moreover, no estimator can do better, in various 
precise mathematical senses (the H\'ajek--LeCam asymptotic local minimax 
theorem and convolution theorem, for instance). Recent work by
\citet{gillmassar00}, surveyed in \citet{gill01}, makes a first 
attempt to carry over these ideas to quantum statistics.

The approach is firstly to delineate more precisely the class of 
attainable
information matrices $I(\theta;M^{(N)})$ based on arbitrary (or special 
classes) of measurements on the model of $N$ 
identical particles each in the same state $\rho(\theta)$. Next, 
using the van Trees inequality,
a Bayesian version of the Cram\'er--Rao 
inequality, see \citet{gilllevit95}, bounds on $I(\theta;M^{(N)})$
are converted into bounds on the asymptotic scaled mean quadratic 
error matrix of regular estimators of $\theta$. Thirdly, one 
constructs measurements and estimators which achieve 
those bounds asymptotically. The first step yields the following theorem.

\begin{thm}[Gill--Massar information bound] In the model of $N$ 
identical copies of a quantum system with state $\rho(\theta)$ on
a $d$ dimensional state space and with $p$ dimensional parameter 
$\theta$,
one has
\begin{equation*}
    \mathrm{trace}(I_{Q}(\theta) ^{-1}I(\theta;M^{(N)})/N)~\le ~
    d-1
\end{equation*}
in any of the following cases: (i) $p=1$ and $d=2$, or (ii) 
$\rho$ is a pure state, or (iii) the measurement $M^{(N)}$ is multi-local.
\end{thm}
A multi-local measurement is a measurement which is composed in an 
arbitrary way of a 
sequence of instruments acting on separate particles. Thus it is 
allowed that the measurement made on particle $2$ depends on the 
outcome of the measurement on particle $1$, and even that after these 
two measurements, yet another measurement, depending on the results 
so far, is made on the first particle in its new state, etc.

In the spin-half case the bound of the above theorem is achievable in the 
sense that for any matrix $J$ such that 
$\mathrm{trace}(I_{Q}(\theta) ^{-1}J)\le 1$, there exists a measurement $M$
on one particle, generally depending on $\theta$, such that
$I(\theta;M)=J$. The measurement is 
a randomised choice of several simple measurements of spin, one 
spin direction for each component of $\theta$.

Application of the van Trees inequality gives the following 
asymptotic bound:

\begin{thm}[Asymptotic information bound] In the model of $N$ 
identical copies of system $\rho(\theta)$,
let $V(\theta)$ denote the limiting scaled mean quadratic error 
matrix of a regular sequence of estimators $\widehat\theta^{(N)}$ 
based on a sequence of measurements $M^{(N)}$ on $N$ particles; i.e.,
$V_{ij}(\theta)=\lim_{N\to\infty}
N \mathrm E_{\theta}\{(\widehat\theta^{(N)}_{i}-\theta_{i})
               (\widehat\theta^{(N)}_{j} - \theta_{j})\}$.
Then $V$ satisfies the inequality
\begin{equation*}
    \mathrm{trace}(I_{Q}(\theta)^{-1} V(\theta)^{-1})~\le~
    d-1
\end{equation*}
in any of the following cases: (i) $p=1$ and $d$=2, or (ii) 
$\rho$ is a pure state, or (iii) the measurements $M^{(N)}$ are 
multi-local.
\end{thm}
A \emph{regular estimator sequence} is one for which the mean quadratic error 
matrices converge uniformly in $\theta$ to a continuous limit. It is
also possible to give a version of the theorem in terms of 
convergence in distribution, H\'ajek-regularity and $V$ the mean quadratic error 
matrix of the limiting distribution, rather than the limit of the 
mean quadratic error.

In the spin-half case, this bound is also asymptically achievable, in 
the sense that for any continuous matrix function $W(\theta)$ such that
$\mathrm{trace}\{I_{Q}(\theta)^{-1} W(\theta)^{-1}\}\le 1$ there exists a sequence of
separable measurements $M^{(N)}$ with asymptotic scaled mean quadratic 
error matrix equal to $W$. This result is proved by consideration of 
a rather natural two-stage measurement procedure. Firstly, on a small
(asymptotically vanishing) proportion of the particles, carry out 
arbitrary measurements allowing consistent estimation of $\theta$,
resulting in a preliminary estimate $\widetilde\theta$.
Then on each of the remaining particles, carry out the measurement 
$\widetilde M$ (on each separate particle) which is optimal in the sense that 
$I(\widetilde\theta;\widetilde M)=J=W(\widetilde\theta)^{-1}$.
Estimate $\theta$ by maximum likelihood estimation, conditional on 
the value of $\widetilde\theta$, on the outcomes obtained in the 
second stage. For large $N$, since $\widetilde\theta$ will then be 
close to the true value of $\theta$, the measurement $\widetilde M$ 
will have Fisher information $I(\theta;\widetilde M)$ close to that of
the `optimal' measurement on one particle with Fisher information
$I(\theta;M)=W(\theta)^{-1}$. By the usual properties of maximum 
likelihood estimators, it will therefore have scaled mean quadratic 
error close to $W(\theta)$.

In the spin-half case we have therefore a complete asymptotic 
efficiency theory in any of the three cases (i) a one-dimensional 
parameter, (ii) a pure state, (iii) multi-local measurements. By 
`complete' we mean that it is precisely known what is the set of all 
attainable limiting scaled mean quadratic error matrices. This 
collection is described in terms of the quantum information matrix for 
one particle.  What is interesting is that when none of these three 
conditions hold, greater asymptotic precision is possible. For 
instance, \citet{gillmassar00} exhibit a generalized measurement 
of two spin-half particles with seven possible outcomes,
which, for a completely unknown mixed state
(a three-parameter model), has about $50\%$ larger total Fisher 
information (for certain parameter values) than any separable 
measurement on two particles. Therefore if one has a large number $N$ 
of particles, one has about $25\%$ better precision when using the 
maximum likelihood estimator applied to the outcomes of this 
measurement on $N/2$ pairs of particles, than any separable 
measurement whatsoever on all $N$. It is not known whether taking 
triples, quadruples, etc., allows even greater increases of precision,
but it seems possible that going to pairs, is enough.
It would be valuable to delineate precisely the set all attainable Fisher 
information matrices when non-separable measurements are allowed on 
each number of particles. 

The measurement in question has seven 
matrix elements. The first is: half the projector onto 
the subspace generated by $\ket{+\vec u_{x}} \otimes \ket{+\vec u_{x}}$.
The next five are obtained by replacing `$+$' with `$-$' and or $x$ 
with $y$ or with $z$. The seventh is the projector onto the state 
spanned by the singlet or Bell state, which we used in teleportation!
This measurement is optimal at the completely mixed state $\rho=\frac12\mathbf 1$ 
for estimating $\vec a$ in the model $\rho=\frac12(\mathbf 1+\vec 
a\cdot\vec{\boldsymbol\sigma})$,
with any loss function which is locally rotation 
invariant.

A similar instance of this phenomenon was called 
\emph{non-locality without entanglement} by \citet{bennettetal99}. One 
could say that though the $N$ particles are not in an entangled state, 
one needs an `entangled measurement', presumably brought about by 
bringing the particles into interaction with one another 
(unitary evolution starting from the product state) before 
measurement, in order to extract maximal information about their 
state. The word `non-locality' refers to the possibility that the $N$ 
particles could be widely separated and brought into interaction 
through other entangled particles; as we saw in Section 3 there are
other examples of this kind in the context of optimal information 
transmission and in teleportation.

\section*{AFTERMATH}

At the beginning of the last (20th) century, two sciences were born and made amazing 
strides: genetics, and quantum physics. In both sciences, randomness plays a central 
part. Whereas this was recognised from the start in genetics---R.A. Fisher is well 
known to biologists, first and foremost, as a great pioneer in genetics---in 
quantum physics this has always been consigned to obscurity, neglect, or 
suppression (Einstein's ``God does not throw dice'').  Now at the beginning 
of the 21st century we seem to be at the threshhold of amazing strides in genetics 
(molecular biology). While biologists claim that this is going to be the century of
molecular biology, physicists argue, also with good reason, that it is going to be 
the century of quantum physics.

I would be surprised if in the coming years, we did not see extraordinary advances in 
physics, in particular, in new quantum technologies. I believe that randomness 
does play a central role in these developments, and that mathematical statisticians 
can and should be involved. I think that the randomness involved in 
quantum mechanics is randomness which can be described by classical probability theory.
However authorities both from physics and from mathematics have 
argued that `quantum probability is a different kind of probability' 
(Feynmann) or that `quantum probability is a strict extension of 
classical probability'. Now in some technical senses this is true. 
If you like to see classical probability theory as part of functional 
analysis, then it is possible to see the algebra of observables of a 
quantum system as a strictly more general kind of mathematical 
structure than the algebra of random variables on a fixed probability 
space. However as we have argued above, the observables of quantum 
mechanics are just an intermediary to deriving the probability 
distributions of outcomes of experiments which could actually be done 
in the laboratory, and then the classical rules of probability theory 
apply. One could just as well argue that the mathematical structure 
of quantum probability theory is a special case of that of classical 
mathematical statistics: it is namely equivalent to a collection of 
classical probability models linked in a rather special way, though 
statisticians are able to consider arbitrary collections of 
probability models.

Another reason why many have claimed that quantum randomness is 
different, is because it can be shown, for instance by considering 
measurements on the two photons in the entangled joint state we used 
for teleportation, that any \emph{deterministic} explanation of the 
randomness in quantum measurements, cannot be described by mere 
statistical variation in hidden or uncontrolled variables of the 
quantum systems, without those variables violating \emph{locality}.
To say it in a different way, any deterministic explanation of what goes 
on in our teleportation example, has to require instantaneous communication from 
Amsterdam to Beijing through the entangled pair of photons.
Though physicists are not happy with randomness, they are even less 
happy with `action at a distance' as this is called.
Thus the randomness in quantum mechanics is of a different nature to 
the randomness in a classical coin toss, for which perfectly 
determinstic rules determine the outcome, starting from the initial 
conditions. However, if one accepts that randomness of a fundamental 
(unexplainable) nature is real, there is not a problem. See my web 
pages for reprints and preprints further discussing these problems.

This brings us to some other philosophical problems, which for me are 
another good reason to be interested in quantum mechanics, since they 
are fundamental issues in physics, deeply connected to probability and 
statistics. We saw in Section 2, three deterministic items concerning 
behaviour of quantum systems `on their own', but a completely 
different and stochastic behaviour when that quantum system is brought 
into interaction with a measurement apparatus, or more generally, with 
the real world.  But a measurement apparatus is just a physical system 
itself, and the system being measured and apparatus doing that, 
together should just be making one large unitary evolution in some 
huge state space. There are no random jumps, no irreversible losses 
of information. Then the same applies if we consider ourselves, the 
observer of the outcome of an experiment, as just another physical 
system\ldots .  This is called the measurement problem. Most working 
physicists are not bothered by it, since they are perfectly able to 
get perfect predictions of experimental results, without worrying 
about the philosophical consistency of the mathematical model. However 
some scientists are deeply worried about it, and there are a number of 
proposals to modify quantum mechanics, or just to modify the interpretation 
rules (the rules by which one draws conclusions from the mathematical 
model back to reality). 

However if one does not attempt to make a mathematical model of all of 
the physical universe, and only models small parts of it, then the 
only problem is a consistency problem, since one might draw the line 
between quantum system and classical outside world, at different levels. 
Working physicists try to draw the border at as macroscopic level as 
possible, and are content that the model prediction at macroscopic 
level is `as if' a quantum jump had taken place to one of several 
macroscopically distinct states, with probabilities which can be 
calculated theoretically, and which are moreover beautifully
confirmed by experiment. It seems to me that careful 
analysis by mathematical statisticians and probabilists could be 
highly valuable, to sift the crazy from the sensible solutions to the 
measurement problem.

\section*{References}
    
\renewcommand{\refname}{}
\ 
\vskip -60pt
\ 

\bibliographystyle{chicago}

\end{document}